\documentclass{iopart}
\usepackage{iopams}
\eqnobysec
\usepackage{graphicx}
\usepackage{epstopdf}
\usepackage{mathrsfs}
\newtheorem{lemma}{Lemma}[section]
\newtheorem{theorem}[lemma]{Theorem}
\newtheorem{remark}{Remark}[section]
\newtheorem{Rule}{Rule}[section]
\newtheorem{proposition}[lemma]{Proposition}

\newtheorem{algorithm}{Algorithm}[section]

\newcommand{\proof} [1]
   { {\bf Proof.} #1 \hfill\opensquare \\}

\def\p{\partial}
\def\d{\delta}

\def\l{\langle}
\def\r{\rangle}

\def\X{\mathcal X}
\def\Y{\mathcal Y}

\def\a{\alpha}

\def\d{\delta}

\def\X{\mathcal X}
\def\Y{\mathcal Y}

\begin{document}
\jl{5}
\title[]
{A fast nonstationary iterative method with convex penalty for inverse problems in Hilbert spaces}

\author{Qinian Jin$^1$ and Xiliang Lu$^2$}

\medskip

\address{$^1$Mathematical Sciences Institute, Australian National University,
Canberra, ACT 0200, Australia}
\address{$^2$School of Mathematics and Statistics, Wuhan University, Wuhan 430072, China}

%\medskip

\eads{Qinian.Jin@anu.edu.au and xllv.math@whu.edu.cn}

%%%%%%%%%%%%%%%%%%%%%%%%%%%%%%%%%%%%%%%%%%%%%%%%%%%%%%%%%%%%%%%%%%%%%%%%%%%%%%%%

%\keywords{Ill-posed problems, inverse problems, regularization,
%Hilbert scales, implicit iteration method, order optimal error
%bounds, general smoothness conditions, operator monotone functions}

\begin{abstract}
  In this paper we consider the computation of approximate solutions for inverse problems in Hilbert spaces. In order to
capture the special feature of solutions, non-smooth convex functions are introduced as penalty terms. By exploiting the
Hilbert space structure of the underlying problems, we propose a fast iterative regularization method which reduces to the classical
nonstationary iterated Tikhonov regularization when the penalty term is chosen to be the square of norm. Each iteration of the method
consists of two steps: the first step involves only the operator from the problem while the second step involves only the penalty term.
This splitting character has the advantage of making the computation efficient. In case the data is corrupted by noise, a stopping
rule is proposed to terminate the method and the corresponding regularization property is established. Finally, we test the performance
of the method by reporting various numerical simulations, including the image deblurring, the determination of source term in Poisson equation,
and the de-autoconvolution problem.
\end{abstract}

%\ams{65J15, 65J20, 47H17}

%\submitto{\IP}

%\maketitle

%\vskip 0.5 cm

%%%%%%%%%%%%%%%%%%%%%% -- SECTION 1 -- %%%%%%%%%%%%%%%%%%%%%%%%%%%%%%%%%%%%%%%%
%%%%%%%%%%%%%%%%%%%%%% -- SECTION 1 -- %%%%%%%%%%%%%%%%%%%%%%%%%%%%%%%%%%%%%%%%
%\section{Introduction}\label{s1}
%

\section{Introduction}

We consider the ill-posed inverse problems of the form
\begin{equation}\label{1.1}
A x=y,
\end{equation}
where $A:\X\to \Y$ is a bounded linear operator between two Hilbert spaces $\X$ and $\Y$ whose inner products and the induced norms
are denoted as $\l \cdot, \cdot\r$ and $\|\cdot\|$ respectively which should be clear from the context. Here the ill-posedness of (\ref{1.1})
refers to the fact that the solution of (\ref{1.1}) does not depend continuously on the data which is a characteristic property of inverse problems.
In practical applications, one never has exact data, instead only noisy data are available due to errors in the measurements. Even if the deviation is very
small, algorithms developed for well-posed problems may fail, since noise could be amplified by an arbitrarily large factor. Therefore,
regularization methods should be used in order to obtain a stable numerical solution. One can refer to \cite{EHN96} for many useful regularization methods
for solving (\ref{1.1}); these methods, however, have the tendency to over-smooth solutions and hence are not quite successful to capture special features.

In case a priori information on the feature of the solution of (\ref{1.1}) is available, we may introduce a proper, lower semi-continuous, convex
function $\Theta: \X \to (-\infty, \infty]$ such that the sought solution of (\ref{1.1}) is in $\mathscr{D}(\Theta)=\{x\in \X: \Theta(x)<\infty\}$.
By taking $x_0\in \X$ and $\xi_0\in \p \Theta(x_0)$, the solution of (\ref{1.1}) with the desired feature can be determined by solving the
constrained minimization problem
\begin{equation}\label{10.14.1}
\min D_{\xi_0} \Theta(x, x_0) \quad \mbox{subject to }  A x =y,
\end{equation}
where $D_{\xi_0}\Theta(x, x_0)$ denotes the Bregman distance induced by $\Theta$ at $x_0$ in the direction $\xi_0$, i.e.
\begin{equation*}
D_{\xi_0} \Theta(x, x_0) = \Theta(x) -\Theta(x_0) -\l \xi_0, x-x_0\r.
\end{equation*}
When only a noisy data $y^\d$ is available, an approximate solution can be constructed by the Tikhonov-type method
\begin{equation}\label{Tik}
x_\a^\d : = \arg \min_{x\in \X} \left\{ \|A x-y^\d\|^2 +\a D_{\xi_0} \Theta(x, x_0)\right\}.
\end{equation}
When the regularization parameter $\a$ is given, many efficient solvers were developed to compute $x_\a^\d$ when $\Theta$ is the
$L^1$ or the total variation function. Unfortunately almost all these methods do not address the choice of $\a$ which, however, is important
for practical applications. In order to use these solvers, one has to perform the trial-and-error procedure to find a reasonable $\a$
which is time consuming. On the other hand, some iterative methods, equipping with proper termination criteria, were developed to find
approximate solutions of (\ref{10.14.1}), see \cite{JZ2012b} and references therein. These iterative methods have the
advantage of avoiding the difficulty for choosing the regularization parameter. However in each iteration step one has to solve a minimization
problem similar to (\ref{Tik}), and overall it may take long time.

In this paper we will propose a fast iterative regularization methods for solving (\ref{10.14.1}) by splitting $A$ and $\Theta$ into different
steps. Our idea is to exploit the Hilbert space structure of the underlying problem to build each iterate by first applying one step of
a well-established classical regularization method and then penalizing the resultant by the convex function $\Theta$. To motivate the method,
we consider the exact data case. We take an invertible bounded linear operator $M: \Y\to \Y$ which can be viewed as a preconditioner
and rewrite (\ref{10.14.1}) into the equivalent form
\begin{equation*}
\min D_{\xi_0}\Theta(x, x_0) \quad \mbox{ subject to } M A x = M y.
\end{equation*}
The corresponding Lagrangian is
\begin{equation*}
{\mathcal L}(x, p) := D_{\xi_0} \Theta(x, x_0) + \l p, M y-M A x\r,
\end{equation*}
where $p\in \Y$ represents the dual variable. Then a desired solution of (\ref{10.14.1}) can be found by determining a saddle point
of ${\mathcal L}$ if exists. Let $(x_c, p_c)$ be a current guess of a saddle point of ${\mathcal L}$, we may update it to get a
new guess $(x_+, p_+)$ as follows:  We first update $p_c$ by solving the proximal maximization problem
\begin{equation*}
p_+:= \arg \max_{p\in \X} \left\{ {\mathcal L}(x_c, p) -\frac{1}{2t} \|p-p_c\|^2\right\}
\end{equation*}
with a suitable step length $t>0$. We then update $x_c$ by solving the minimization problem
\begin{equation*}
x_+ := \arg \min_{x\in \X} {\mathcal L}(x, p_+).
\end{equation*}
By straightforward calculation and simplification it follows
\begin{eqnarray*}
p_+ &= p_c - t M (A x_c-y),\\
x_+ &= \arg \min_{x\in \X} \left\{ \Theta(x) -\l \xi_0+A^* M^* p_+, x\r\right\}
\end{eqnarray*}
which is the one step result of the Uzawa algorithm \cite{AHU58} or the dual subgradient method \cite{Shor85}, where $A^*:\Y \to \X$
and $M^*: \Y \to \Y$ denote the adjoint operators of $A$ and $M$ respectively. By setting $\xi_c = \xi_0+A^* M^* p_c$ and $\xi_+ = \xi_0 + A^*M^* p_+$,
the above equation can be transformed into the form
\begin{equation}\label{10.13.2}
\left\{\begin{array}{lll}
\xi_+ = \xi_c - t A^*M^* M (A x_c-y),\\
x_+ = \arg \min_{x\in \X} \left\{ \Theta(x) -\l \xi_+, x\r\right\}.
\end{array}\right.
\end{equation}

Now we may apply the updating scheme (\ref{10.13.2}) to $A x =y$ iteratively but with dynamic preconditioning operator $M_n:\Y\to \Y$
and variable step size $t_n>0$. This gives rise to the following iterative methods
\begin{eqnarray}\label{method0}
\left\{\begin{array}{lll}
\xi_{n+1}=\xi_n - t_n A^*M_n^* M_n (A x_n-y),\\
x_{n+1} = \arg \min_{x\in \X} \left\{ \Theta(x) -\l \xi_{n+1}, x\r \right\}.
\end{array}\right.
\end{eqnarray}
The performance of the method (\ref{method0}) depends on the choices of $\{M_n\}$. If we take $M_n =I$
for all $n$, (\ref{method0}) becomes the method that has been studied in \cite{BH2012,JW2013} which is the generalization of the classical
Landweber iteration and is known to be a slowly convergent method.

In this paper we will consider (\ref{method}) with $M_n =(\a_n I +AA^*)^{-1/2}$ for all $n$, where $\{\a_n\}$ is a decreasing
sequence of positive numbers. This yields the nonstationary iterative method
\begin{eqnarray}\label{method}
\left\{\begin{array}{lll}
\xi_{n+1}=\xi_n - t_n A^*(\a_n I +AA^*)^{-1} (A x_n-y),\\
x_{n+1} = \arg \min_{x\in \X} \left\{ \Theta(x) -\l \xi_{n+1}, x\r \right\}.
\end{array}\right.
\end{eqnarray}
Observing that when $\Theta(x) =\|x\|^2/2$ and $t_n =1$ for all $n$, (\ref{method}) reduces to the nonstationary iterated Tikhonov
regularization
\begin{equation}\label{NSIT}
x_{n+1} =x_n - A^* (\a_n I +AA^*)^{-1} (A x_n -y)
\end{equation}
whose convergence has been studied in \cite{HG98} and it has been shown to be a fast convergent method when $\{\a_n\}$ is a
geometric decreasing sequence. This strongly suggests that our method (\ref{method}) may also exhibit fast convergence property if
$\{\a_n\}$ and $\{t_n\}$ are chosen properly. We will confirm this in the present paper. It is worthy to point out that each iteration
in (\ref{method}) consists of two steps: the first step involves only the operator $A$ and the second step involves only the convex
function $\Theta$. This splitting character can make the computation much easier.

This paper is organized as follows. In section 2, we start with some preliminary facts from convex analysis, and then give the convergence
analysis of the method (\ref{method}) when the data is given exactly. In case the data is corrupted by noise, we propose a stopping
rule to terminate the iteration and establish the regularization property. We also give a possible extension of our method to solve
nonlinear inverse problems in Hilbert spaces. In section 3 we test the performance of our method by reporting various numerical simulations,
including the image deblurring, the determination of source term in Poisson equation and the de-autoconvolution problem.

\section{Convergence analysis of the method}

In this section we first give the convergence analysis of (\ref{method}) with suitable chosen $t_n$ when $\Theta: \X \to (-\infty, \infty]$ is
a proper, lower semi-continuous function that is strongly convex in the sense that there is a constant $c_0>0$ such that
\begin{equation}\label{sc}
\Theta(s \bar x +(1-s) x)+c_0 s(1-s) \|\bar x-x\|^2 \le s \Theta(\bar x) +(1-s) \Theta(x)
\end{equation}
for all $0\le s\le 1$ and $\bar x, x\in \X$. We then consider the method when the data contains noise and propose a stopping rule to render
it into a regularization method. Our analysis is based on some important results from convex analysis which will be recalled in the
following subsection.

\subsection{\bf Tools from convex analysis}

Given a convex function $\Theta: \X \to (-\infty, \infty]$, we will use $\mathscr{D}(\Theta): =\{x\in \X, \Theta(x) <\infty\}$
to denote its effective domain. It is called proper if $\mathscr{D}(\Theta)\ne \emptyset$. Given $x\in \X$, the set
\begin{equation*}
\p \Theta(x) :=\{\xi \in \X: \Theta(\bar x)-\Theta(x) -\l \xi, \bar x- x\r\ge 0  \mbox{ for all } \bar x \in \X\}
\end{equation*}
is called the subdifferential of $\Theta$ at $x$ and each element $\xi\in \p \Theta(x)$ is called a subgradient.

Our convergence analysis of (\ref{method}) will not be carried out directly under the norm of $\X$. Instead we will use the Bregman distance
(\cite{Br1967}) induced by $\Theta$. Given $\xi\in \p \Theta(x)$, the quantity
\begin{equation*}
D_{\xi}\Theta(\bar x, x) := \Theta(\bar x) -\Theta(x) -\l \xi, \bar x-x\r, \quad \bar x \in \X
\end{equation*}
is called the Bregman distance induced by $\Theta$ at $x$ in the direction $\xi$. It is clear that $D_\xi \Theta(\bar x, x) \ge 0$.
However, Bregman distance is not a metric distance since it does not satisfy the symmetry and the triangle inequality in general. Nevertheless, when
$\Theta$ is strongly convex in the sense of (\ref{sc}), there holds (\cite{Z2002})
\begin{equation*}
D_{\xi}\Theta(\bar x, x) \ge c_0 \|\bar x-x\|^2, \qquad \forall \bar x\in \X \mbox{ and } \xi \in \p \Theta(x)
\end{equation*}
which means that the Bregman distance can be used to detect information under the norm of $\X$.

Although $\Theta$ could be non-smooth, its Fenchel conjugate can have enough regularity provided $\Theta$ has enough convexity.
The Fenchel conjugate of $\Theta$ is defined by
\begin{equation*}
\Theta^*(\xi) := \sup_{x\in \X} \left\{ \l \xi, x\r -\Theta(x)\right\}, \qquad \forall \xi\in \X.
\end{equation*}
For a proper, lower semi-continuous, convex function $\Theta$, there always holds
\begin{equation*}
\xi\in \p \Theta(x) \Longleftrightarrow x\in \p \Theta^*(\xi) \Longleftrightarrow \Theta(x) +\Theta^*(\xi) =\l \xi, x\r.
\end{equation*}
Consequently, the Bregman distance can be equivalently written as
\begin{equation}\label{10.13.3}
D_{\xi}\Theta(\bar x, x) = \Theta(\bar x) +\Theta^*(\xi) -\l \xi, \bar x\r.
\end{equation}
If in addition $\Theta$ is strongly convex in the sense of (\ref{sc}), then $\mathscr{D}(\Theta^*)=\X$, $\Theta^*$ is Fr\'{e}chet differentiable,
and its gradient $\nabla \Theta^*$ satisfies
\begin{equation}\label{10.13.4}
\|\nabla \Theta^*(\xi)-\nabla \Theta^*(\eta)\|\le \frac{\|\xi-\eta\|}{2c_0},
\end{equation}
i.e. $\nabla \Theta^*$ is Lipschitz continuous. These facts are crucial in the forthcoming convergence analysis and their
proofs can be found in many standard textbooks, cf. \cite{Z2002}.

\subsection{\bf The method with exact data}

We consider the convergence of the method (\ref{method}) under the condition that $\Theta$ is proper, lower semi-continuous, and strongly convex
in the sense of (\ref{sc}). We will always assume that (\ref{1.1}) has a solution in $\mathscr{D}(\Theta)$. By taking $\xi_0\in \X$ and define
\begin{equation*}
x_0 =\arg \min_{x\in \X} \left\{ \Theta(x) -\l \xi_0, x\r \right\}
\end{equation*}
as an initial guess, we define $x^\dag$ to be the solution of (\ref{1.1}) in $\mathscr{D}(\Theta)$ satisfying
\begin{equation}\label{10.13.1}
D_{\xi_0}\Theta(x^\dag, x_0) = \min\left\{ D_{\xi_0} \Theta(x, x_0): A x =y\right\}.
\end{equation}
It is easy to show that such $x^\dag$ is uniquely defined. Our aim is to show that the sequence $\{x_n\}$ produced by (\ref{method}) eventually
converges to $x^\dag$ if $t_n$ is chosen properly.

To this end, we first consider the monotonicity of the Bregman distance $D_{\xi_n}\Theta(\hat x, x_n)$ with respect to $n$ for any
solution $\hat x$ of (\ref{1.1}) in $\mathscr{D}(\Theta)$. By the subdifferential calculus and the definition of $x_n$, it is easy to see that
$\xi_n \in \p \Theta(x_n)$ and hence $x_n =\nabla \Theta^*(\xi_n)$. Therefore, in view of (\ref{10.13.3}) and (\ref{10.13.4}) we have
\begin{eqnarray*}
\fl D_{\xi_{n+1}}\Theta(\hat x, x_{n+1}) -D_{\xi_n}\Theta(\hat x, x_n)
 = \Theta^*(\xi_{n+1}) -\Theta^*(\xi_n) -\l \xi_{n+1}-\xi_n, \hat x\r\\
 = \Theta^*(\xi_{n+1}) -\Theta^*(\xi_n) -\l \xi_{n+1}-\xi_n, \nabla \Theta^*(\xi_n)\r -\l \xi_{n+1}-\xi_n, \hat x-x_n\r\\
 =\int_0^1 \l \xi_{n+1}-\xi_n, \nabla \Theta^*(\xi_n+s(\xi_{n+1}-\xi_n))-\nabla \Theta^*(\xi_n)\r ds\\
 \quad \,  -\l \xi_{n+1}-\xi_n, \hat x-x_n\r\\
 \le \int_0^1 \frac{1}{2 c_0} s \|\xi_{n+1}-\xi_n\|^2 ds -\l \xi_{n+1}-\xi_n, \hat x-x_n\r\\
 = \frac{1}{4 c_0} \|\xi_{n+1}-\xi_n\|^2 -\l \xi_{n+1}-\xi_n, \hat x -x_n\r.
\end{eqnarray*}
Using the definition of $\xi_{n+1}$ in (\ref{method}) and $A \hat x =y$ we obtain
\begin{eqnarray*}
\fl D_{\xi_{n+1}}\Theta(\hat x, x_{n+1}) -D_{\xi_n}\Theta(\hat x, x_n)
& \le \frac{1}{4 c_0} t_n^2 \|A^*(\a_n I + AA^*)^{-1} (A x_n-y)\|^2 \\
& \quad\, -t_n \l (\a_n I + AA^*)^{-1} (A x_n-y), A x_n-y\r.
\end{eqnarray*}
If $A x_n -y \ne 0$, we may choose $t_n$ such that
\begin{equation}\label{10.13.5}
t_n = \frac{\mu_0 \l(\a_n I + AA^*)^{-1}(A x_n-y), A x_n -y\r}{\|A^*(\a_n I + AA^*)^{-1}(A x_n-y)\|^2}
\end{equation}
with $0<\mu_0<4 c_0$, then it yields
\begin{eqnarray} \label{eq:11.22}
&D_{\xi_{n+1}}  \Theta(\hat x, x_{n+1}) -D_{\xi_n}\Theta(\hat x, x_n) \nonumber\\
&\le -\left(1-\frac{\mu_0}{4 c_0}\right) t_n \| (\a_n I + AA^*)^{-1/2} (A x_n-y)\|^2 \le 0.
\end{eqnarray}
When $A x_n -y=0$, the inequality (\ref{eq:11.22}) obviously holds for any $t_n\ge 0$.  We observe that the $t_n$ chosen by
(\ref{10.13.5}) could be very large when $\|A x_n-y\|$ is small. Using such a choice of $t_n$ it could make the method
numerically unstable, in particular when the data contains noise. To avoid this, we take a preassigned
number $\mu_1>0$ and then set
\begin{equation}\label{10.13.7}
t_n = \min \left\{\frac{\mu_0\l (\a_n I + AA^*)^{-1} (A x_n-y), A x_n -y\r}{\|A^*(\a_n I + AA^*)^{-1} (A x_n-y)\|^2}, \mu_1\right\}.
\end{equation}
The above argument then shows the following monotonicity result.

\begin{lemma}\label{lem10.13.1}
If $t_n$ is chosen by (\ref{10.13.7}) with $0<\mu_0<4 c_0$ and $\mu_1>0$, then
\begin{equation*}
D_{\xi_{n+1}}  \Theta(\hat x, x_{n+1}) \le D_{\xi_n}\Theta(\hat x, x_n)
\end{equation*}
and
\begin{equation}\label{10.3.4}
\fl \qquad c_1 t_n \|(\a_n I + AA^*)^{-1/2} (A x_n-y)\|^2 \le D_{\xi_n} \Theta(\hat x, x_n) -D_{\xi_{n+1}}\Theta(\hat x, x_{n+1})
\end{equation}
for any solution $\hat x$ of (\ref{1.1}) in $\mathscr{D}(\Theta)$, where $c_1:= 1-\mu_0/(4 c_0)$.
\end{lemma}

We will use Lemma \ref{lem10.13.1} to derive the convergence of the method (\ref{method}).  For the step size $t_n$ defined by (\ref{10.13.7}),
it is easy to see that
\begin{equation*}
\min\{\mu_0, \mu_1\} \le t_n \le \mu_1,
\end{equation*}
where we used the inequality $\|A^*(\a_n I +AA^*)^{-1/2}\|\le 1$ to derive the left inequality. This together with (\ref{10.3.4}) implies
\begin{equation}\label{10.13.8}
\lim_{n \rightarrow \infty} \|(\a_n I + AA^*)^{-1/2} (A x_n-y)\| =0.
\end{equation}
Since $\|(\a_n I + AA^*)^{1/2}\|\le \sqrt{\a_0+\|A\|^2}$, we can further derive that
\begin{equation}\label{10.13.9}
\lim_{n \rightarrow \infty} \|A x_n-y\| =0.
\end{equation}

The following main result shows that the method (\ref{method}) is indeed convergent if $t_n$ is chosen by (\ref{10.13.7}).

\begin{theorem}\label{thm1}
Let $\Theta: \X \to (-\infty, \infty]$ be a a proper, lower semi-continuous function that is strongly convex in the sense of (\ref{sc}).
If $\{\a_n\}$ is a decreasing sequence of positive numbers and if $t_n$ is chosen by (\ref{10.13.7}) with $0<\mu_0<4 c_0$ and $\mu_1>0$, then for the
method (\ref{method}) there hold
\begin{equation*}
\lim_{n\rightarrow \infty} \|x_n-x^\dag\|=0 \quad \mbox{and} \quad \lim_{n\rightarrow \infty} D_{\xi_n} \Theta(x^\dag, x_n)=0.
\end{equation*}

\end{theorem}

The proof is based on the following useful result.

\begin{proposition} \label{general}
Consider the equation (\ref{1.1}). Let $\Theta: \X\to (-\infty, \infty]$ be a proper, lower semi-continuous and strong convex function.
Let $\{x_n\}\subset \X$ and $\{\xi_n\}\subset \X$ be such that

\begin{enumerate}

\item[(i)] $\xi_n\in \p \Theta(x_n)$ for all $n$;

\item[(ii)] for any solution $\hat x$ of (\ref{1.1}) in $\mathscr{D}(\Theta)$ the sequence $\{D_{\xi_n}\Theta(\hat x, x_n)\}$ is monotonically decreasing;

\item[(iii)] $\lim_{n\rightarrow \infty} \|A x_n-y\|=0$.

\item[(iv)] there is a subsequence $\{n_k\}$ with $n_k\rightarrow \infty$ such that for any solution $\hat x$ of (\ref{1.1})
in $\mathscr{D}(\Theta)$ there holds
\begin{equation}\label{10.3.2}
 \lim_{l\rightarrow \infty} \sup_{k\ge l} |\l \xi_{n_k}-\xi_{n_l}, x_{n_k}-\hat x\r| =0.
\end{equation}
\end{enumerate}
Then there exists a solution $x_*$ of (\ref{1.1}) in $\mathscr{D}(\Theta)$ such that
\begin{equation*}
\lim_{n\rightarrow \infty} D_{\xi_n}\Theta(x_*, x_n)=0.
\end{equation*}
If, in addition, $\xi_{n+1}-\xi_n \in \mathscr{R}(A^*)$ for all $n$, then $x_*=x^\dag$.
\end{proposition}

\proof{
This is a slight modification of \cite[Proposition 3.6]{JW2013}, we include here the proof for completeness.

We first show the convergence of $\{x_{n_k}\}$. For any $l<k$ we have from the definition of Bregman distance that
\begin{equation}\label{10.3.1}
\fl D_{\xi_{n_l}}\Theta(x_{n_k}, x_{n_l})  =D_{\xi_{n_l}}\Theta(\hat{x}, x_{n_l}) -D_{\xi_{n_k}}\Theta(\hat{x}, x_{n_k})
+\l \xi_{n_k}-\xi_{n_l}, x_{n_k}-\hat{x}\r.
\end{equation}
By the monotonicity of $\{D_{\xi_n}\Theta(\hat x, x_n)\}$  and (\ref{10.3.2}) we obtain that $D_{\xi_{n_l}}\Theta(x_{n_k}, x_{n_l})\rightarrow 0$
as $k, l\rightarrow \infty$. In view of the strong convexity of $\Theta$, it follows that $\{x_{n_k}\}$ is a Cauchy sequence in $\X$.
Thus $x_{n_k}\rightarrow x_*$ for some $x_*\in \X$. Since $\lim_{n\rightarrow \infty} \|A x_n -y\|=0$, we have $A x_*=y$.

In order to show $x_*\in \mathscr{D}(\Theta)$, we use $\xi_{n_k}\in \p \Theta(x_{n_k})$ to obtain
\begin{equation}\label{5.5.3.1}
\Theta(x_{n_k})\le \Theta(\hat{x}) +\l \xi_{n_k}, x_{n_k}-\hat{x}\r.
\end{equation}
In view of (\ref{10.3.2}) and $x_{n_k}\rightarrow x_*$ as $k\rightarrow \infty$, there is a constant $C_0$ such that
\begin{equation*}
|\l \xi_{n_k} -\xi_{n_0}, x_{n_k}-\hat x\r| \le C_0 \quad \mbox{and} \quad |\l \xi_{n_0}, x_{n_k}-\hat x\r| \le C_0, \quad \forall k.
\end{equation*}
Therefore $|\l \xi_{n_k}, x_{n_k}-\hat{x}\r|\le 2C_0$ for all $k$. By using the lower semi-continuity of $\Theta$ we obtain
\begin{equation*}
\Theta(x_*)\le \liminf_{k\rightarrow\infty} \Theta(x_{n_k})\le
\Theta(\hat{x}) + 2C_0 <\infty.
\end{equation*}
This implies that $x_*\in \mathscr{D}(\Theta)$.

Next we derive the convergence in Bregman distance.  Since $\{D_{\xi_n} \Theta(x_*, x_n)\}$ is monotonically decreasing, the limit
$\varepsilon_0:=\lim_{n\rightarrow \infty} D_{\xi_n}\Theta(x_*, x_n)\ge 0$ exists. By taking $k\rightarrow \infty$ in (\ref{10.3.1}) with $\hat x=x_*$
and using the lower semi-continuous of $\Theta$, we obtain
\begin{equation*}
D_{\xi_{n_l}} \Theta(x_*, x_{n_l}) \le D_{\xi_{n_l}}\Theta(x_*, x_{n_l}) -\varepsilon_0 + \sup_{k\ge l} |\l \xi_{n_k}-\xi_{n_l}, x_{n_k}- x_*\r|
\end{equation*}
which is true for all $l$. Letting $l\rightarrow \infty$ and using (\ref{10.3.2}) gives $\varepsilon_0 \le \varepsilon_0-\varepsilon_0=0$.
Thus $\varepsilon_0=0$, i.e. $\lim_{n\rightarrow \infty} D_{\xi_n} \Theta(x_*, x_n)=0$.

Finally we show that $x_*=x^\dag$. We use (\ref{5.5.3.1}) with $\hat{x}$ replaced by $x^\dag$ to obtain
\begin{equation}\label{5.5.7}
D_{\xi_0}\Theta(x_{n_k}, x_0)\le D_{\xi_0}\Theta(x^\dag, x_0)+\l \xi_{n_k}-\xi_0, x_{n_k}-x^\dag\r.
\end{equation}
Because of (\ref{10.3.2}), for any $\varepsilon>0$  we can find $k_0$ such that
\begin{equation*}
\left|\l \xi_{n_k}-\xi_{n_{k_0}}, x_{n_k}-x^\dag\r\right| <\varepsilon/2, \qquad \forall k\ge k_0.
\end{equation*}
We next consider $\l \xi_{n_{k_0}}-\xi_0, x_{n_k}-x^\dag\r$. Since $\xi_{n+1}-\xi_n \in \mathscr{R}(A^*)$,
we can find $v\in \Y$ such that $\xi_{n_{k_0}}-\xi_0=A^* v$. Consequently
\begin{equation*}
|\l \xi_{n_{k_0}}-\xi_0, x_{n_k}-x^\dag\r | =|\l v, A x_{n_k}-y\r| \le \|v\| \|A x_{n_k}-y\|.
\end{equation*}
Since $\|A x_n-y\|\rightarrow 0$ as $n\rightarrow \infty$, we can find $k_1\ge k_0$ such that
\begin{equation*}
|\l \xi_{n_{k_0}}-\xi_0, x_{n_k}-x^\dag\r |<\varepsilon/2, \qquad \forall k\ge k_1.
\end{equation*}
Therefore $|\l \xi_{n_k}-\xi_0, x_{n_k}-x^\dag\r|<\varepsilon$ for all $k\ge k_1$. Since $\varepsilon>0$ is arbitrary,
we obtain $\lim_{k\rightarrow \infty} \l \xi_{n_k}-\xi_0, x_{n_k}-x^\dag\r=0$. By taking $k\rightarrow \infty$ in
(\ref{5.5.7}) and using the lower semi-continuity of $\Theta$ we obtain $D_{\xi_0}\Theta(x_*, x_0)\le D_{\xi_0}\Theta(x^\dag, x_0)$.
According to the definition of $x^\dag$ we must have $D_{\xi_0}\Theta(x_*, x_0)=D_{\xi_0}\Theta(x^\dag, x_0)$.
By uniqueness it follows $x_*=x^\dag$.
}

\noindent{\bf Proof of Theorem \ref{thm1}.}
We will use Proposition \ref{general} to complete the proof. By the definition of $\{\xi_n\}$ we always have $\xi_{n+1}-\xi_n\in \mathscr{R}(A^*)$.
It remains to verify the four conditions in Proposition \ref{general}. By the definition of $x_n$ we have $\xi_n \in \p \Theta(x_n)$
which implies (i) in Proposition \ref{general}. Moreover, Lemma \ref{lem10.13.1} and (\ref{10.13.9}) confirm (ii) and (iii) in Proposition \ref{general}
respectively.

It remains only to verify (iv) in Proposition \ref{general}. To this end, we consider
\begin{equation*}
R_n : = \| (\a_n I + AA^*)^{-1/2} (A x_n-y)\|.
\end{equation*}
In view of (\ref{10.13.8}), we have $\lim_{n\rightarrow \infty} R_n=0$. Moreover, by the definition of the method (\ref{method}),
if $R_n=0$ for some $n$, then $R_m=0$ for all $m\ge n$. Consequently, we may choose a strictly increasing subsequence $\{n_k\}$ of integers such that
$n_0=0$ and $n_k$, for each $k\ge 1$, is the first integer satisfying
\begin{equation*}
n_k \ge n_{k-1}+1  \quad \mbox{ and } \quad R_{n_k} \le R_{n_{k-1}}.
\end{equation*}
For this chosen $\{n_k\}$ it is easy to see that
\begin{equation}\label{9.29.1}
R_n \ge R_{n_k}, \qquad \forall 0\le n\le n_k.
\end{equation}
Inddeed, for $0\le n<n_k$, we can find $0\le l< k$ such that $n_l \le n < n_{l+1}$ and thus,
by the definition of $n_{l+1}$, we have $R_n\ge R_{n_l} \ge R_{n_k}$.
With the above chosen $\{n_k\}$, we will show that (\ref{10.3.2}) holds for any solution $\hat x$ of (\ref{1.1}) in $\mathscr{D}(\Theta)$.
By the definition of $\xi_n$ we have for $l<k$ that
\begin{eqnarray*}
\l \xi_{n_l}-\xi_{n_k}, \hat x-x_{n_k}\r &=\sum_{n=n_l}^{n_k-1} \l \xi_{n+1}-\xi_n, x_{n_k}-\hat x\r\\
& = - \sum_{n=n_l}^{n_k-1} t_n \l (\a_n I + AA^*)^{-1} (A x_n-y), A x_{n_k}-y\r.
\end{eqnarray*}
Therefore
\begin{eqnarray*}
&|\l \xi_{n_l}-\xi_{n_k}, \hat x-x_{n_k}\r| \\
&\le \sum_{n=n_l}^{n_k-1} t_n \| (\a_n I + AA^*)^{-1/2} (A x_n-y)\| \|(\a_n I + AA^*)^{-1/2} (A x_{n_k}-y)\|.
\end{eqnarray*}
By using the monotonicity of $\{\a_n\}$ and (\ref{9.29.1}), we have for $0\le n\le n_k$ that
\begin{equation*}
\fl \|(\a_n I + AA^*)^{-1/2} (A x_{n_k}-y)\| \le \|(\a_{n_k} I + AA^*)^{-1/2} (A x_{n_k}-y)\| =R_{n_k} \le R_n.
\end{equation*}
Consequently, it follows from  (\ref{10.3.4}) that
\begin{eqnarray*}
\fl |\l \xi_{n_l}-\xi_{n_k}, \hat x-x_{n_k}\r| &\le \sum_{n=n_l}^{n_k-1} t_n R_n^2
\le \frac{1}{c_1}\left( D_{\xi_{n_l}}\Theta(\hat x, x_{n_l}) -D_{\xi_{n_k}} \Theta(\hat x, x_{n_k})\right)
\end{eqnarray*}
which, together with the monotonicity of $\{D_{\xi_n}\Theta(\hat x, x_n)\}$, implies (\ref{10.3.2}). The proof is therefore complete. \hfill $\Box$

\subsection{\bf The method with noisy data}

We next consider the situation that the data contains noise. Thus, instead of $y$, we only have noisy data $y^\d$ satisfying
\begin{equation*}
\|y^\d-y\| \le \d
\end{equation*}
with a small known noise level $\d>0$. The corresponding method takes the form
\begin{equation}\label{method_noise}
\left\{\begin{array}{lll}
\xi_{n+1}^\d =\xi_n^\d - t_n^\d A^*(\a_n I + AA^*)^{-1} (A x_n^\d-y^\d),\\
x_{n+1}^\d = \arg \min_{x\in \X} \left\{ \Theta(x) -\l \xi_{n+1}^\d, x\r \right\}
\end{array}\right.
\end{equation}
with suitably chosen step length $t_n^\d>0$, where $\xi_0^\d:=\xi_0$ and $x_0^\d :=x_0$. In order to terminate the method, we need some stopping criterion.
It seems that a natural one is the discrepancy principle
\begin{equation}\label{DP}
\|A x_{n_\d}^\d -y^\d\| \le \tau \d <\|A x_n^\d -y^\d\|, \qquad 0\le n< n_\d
\end{equation}
for some number $\tau>1$. Unfortunately, we can not prove the regularization property for the method terminated by the discrepancy principle;
furthermore, numerical simulations indicate that the discrepancy principle might not always produce satisfactory reconstruction result.
Therefore, the discrepancy principle might not be a natural rule to terminate (\ref{method_noise}). Recall that when we motivate our method,
we consider the preconditioned equation
\begin{equation*}
(\a_n I + AA^*)^{-1/2} A x =(\a_n I + AA^*)^{-1/2} y
\end{equation*}
instead of $A x =y$. This indicates that it might be natural to stop the iteration as long as
\begin{equation}\label{PMD}
\|(\a_n I + AA^*)^{-1/2} (A x_n^\d -y^\d)\| \le \tau \|(\a_n I + AA^*)^{-1/2}(y-y^\d)\|
\end{equation}
is satisfied for the first time. The stopping rule (\ref{PMD}) can be viewed as the discrepancy principle applied to the preconditioned equation.
Since the right hand side of (\ref{PMD}) involves $y$ which is not available, it can not be used in practical applications.
Considering $\|(\a_n I +AA^*)^{-1/2}\|\le 1/\sqrt{\a_n}$, we may replace the right hand side of (\ref{PMD}) by $\tau \d/\sqrt{\a_n}$
which leads to the following stopping rule.

\begin{Rule}\label{Rule1}
Let $\tau>1$ be a given number. We define $n_\d$ to be the first integer such that
\begin{equation*}
\a_{n_\d} \l (\a_{n_\d} I + AA^*)^{-1} (A x_{n_\d}^\d-y^\d), A x_{n_\d}^\d-y^\d\r \le \tau^2 \d^2.
\end{equation*}
\end{Rule}

In the context of Tikhonov regularization for linear ill-posed inverse problems, a similar rule was proposed in \cite{Rau84,Gfr87}
to choose the regularization parameter. The rule was then generalized and analyzed in \cite{SEK93,JH99} for nonlinear Tikhonov regularization
and was further extended in \cite{Jin2000} as a stopping rule for the iteratively regularized Gauss-Newton method
for solving nonlinear inverse problems in Hilbert spaces.

%It is easy to see that Rule \ref{Rule1} always terminates the
%iteration (\ref{method_noise}) no later than the discrepancy principle (\ref{DP}).

Combining Rule \ref{Rule1} with (\ref{method_noise}) and using suitable choice of the step length $t_n^\d$ it yields the following algorithm.

\begin{algorithm} [Nonstationary iterative method with convex penalty] \label{alg1} \quad

\begin{enumerate}
\item[(i)] Take $\tau>1$, $\mu_0>0$, $\mu_1>0$ and a decreasing sequence $\{\a_n\}$ of positive numbers;

\item[(ii)] Take $\xi_0\in \X$ and define $x_0:=\arg\min_{x\in \X} \{\Theta(x) -\l \xi_0, x\r\}$ as an initial guess;

\item[(iii)]  For each $n=0, 1, \cdots$ define $\xi_{n+1}^\d$ and $x_{n+1}^\d$ by (\ref{method_noise}), where
\begin{equation}\label{10.3.6}
t_n^\d = \min \left\{\frac{\mu_0\l (\a_n I + AA^*)^{-1}(A x_n^\d-y^\d), A x_n^\d-y^\d\r}{\|A^*(\a_n I + AA^*)^{-1} (A x_n^\d-y^\d)\|^2}, \mu_1\right\}
\end{equation}

\item[(iv)] Let $n_\d$ be the integer determined by Rule \ref{Rule1} and use $x_{n_\d}^\d$ as an approximate solution.
\end{enumerate}
\end{algorithm}

 The following lemma shows that Algorithm \ref{alg1} is well defined and certain monotonicity result holds along the iteration
 if $\mu_0>0$ is suitably small.

\begin{lemma}\label{lem1}
Let $\Theta: \X \to (-\infty, \infty]$ be a a proper, lower semi-continuous function that is strongly convex in the sense of (\ref{sc}).
If $\{\a_n\}$ is a decreasing sequence of positive numbers and $t_n^\d$ is chosen by (\ref{10.3.6}) with $0<\mu_0<4 c_0 (1-1/\tau)$
and $\mu_1>0$, then Rule \ref{Rule1} defines a finite integer $n_\d$. Moreover, if $n_\d\ge 1$, then for the sequences $\{\xi_n^\d\}$ and $\{x_n^\d\}$
defined by (\ref{method_noise}) there holds
\begin{equation}\label{monotone}
D_{\xi_{n+1}^\d}\Theta (\hat x, x_{n+1}^\d) \le D_{\xi_n^\d} \Theta (\hat x, x_n^\d), \qquad 0\le n<n_\d
\end{equation}
for any solution $\hat x$ of (\ref{1.1}) in $\mathscr{D}(\Theta)$.
\end{lemma}

\proof{
Let $0\le n<n_\d$. By using the similar argument in the proof of Lemma \ref{lem10.13.1} we can obtain
\begin{eqnarray*}
\fl D_{\xi_{n+1}^\d}\Theta(\hat x, x_{n+1}^\d) -D_{\xi_n^\d}\Theta(\hat x, x_n^\d)
& \le \frac{1}{4 c_0} \|\xi_{n+1}^\d-\xi_n^\d\|^2 -\l \xi_{n+1}^\d-\xi_n^\d, \hat x -x_n^\d\r\\
& = \frac{1}{4 c_0} (t_n^\d)^2 \|A^*(\a_n I + AA^*)^{-1} (A x_n^\d-y^\d)\|^2\\
& \quad \, -t_n^\d \l (\a_n I + AA^*)^{-1} (A x_n^\d-y^\d), A x_n^\d -y\r.
\end{eqnarray*}
In view of $\|y^\d-y\|\le \d$ and the choice of $t_n^\d$, it follows that
\begin{eqnarray*}
\fl D_{\xi_{n+1}^\d}  \Theta(\hat x, x_{n+1}^\d) -D_{\xi_n^\d}\Theta(\hat x, x_n^\d)
& \le -\left(1-\frac{\mu_0}{4 c_0}\right) t_n^\d \| (\a_n I + AA^*)^{-1/2} (A x_n^\d-y^\d)\|^2 \\
& \quad \, +  t_n^\d \| (\a_n I + AA^*)^{-1/2}(A x_n^\d-y^\d)\| \frac{\d}{\sqrt{\a_n}}.
\end{eqnarray*}
By the definition of $n_\d$ and $n<n_\d$ we have
\begin{equation}\label{10.3.5}
\frac{\d}{\sqrt{\a_n}} \le \frac{1}{\tau} \|(\a_n I + AA^*)^{-1/2} (A x_n^\d -y^\d)\|.
\end{equation}
Therefore, we have with $c_2 := 1-1/\tau -\mu_0/(4 c_0)>0$ that
\begin{eqnarray*}
\fl \qquad D_{\xi_{n+1}^\d}  \Theta(\hat x, x_{n+1}^\d) -D_{\xi_n^\d}\Theta(\hat x, x_n^\d)
\le -c_2 t_n^\d \|(\a_n I + AA^*)^{-1/2} (A x_n^\d-y^\d)\|^2 \le 0.
\end{eqnarray*}
This shows the monotonicity result (\ref{monotone}) and
\begin{equation*}
\fl \qquad c_2 t_n^\d \|(\a_n I + AA^*)^{-1/2} (A x_n^\d-y^\d)\|^2 \le D_{\xi_n^\d} \Theta(\hat x, x_n^\d) -D_{\xi_{n+1}^\d}\Theta(\hat x, x_{n+1}^\d)
\end{equation*}
for all $0\le n<n_\d$. We may sum the above inequality over $n$ from $n=0$ to $n=m$ for any $m<n_\d$ to get
\begin{equation*}
c_2 \sum_{n=0}^m t_n^\d \|(\a_n I + AA^*)^{-1/2} (A x_n^\d-y^\d)\|^2 \le D_{\xi_0} \Theta(\hat x, x_0).
\end{equation*}
By the choice of $t_n^\d$ it is easy to check that $t_n^\d\ge \min\{\mu_0, \mu_1\}$.
Therefore, in view of (\ref{10.3.5}), we have
\begin{equation}\label{10.6.2}
c_2 \min\{\mu_0, \mu_1\} \tau^2 \d^2 \sum_{n=0}^m \frac{1}{\a_n} \le D_{\xi_0}\Theta(\hat x, x_0)
\end{equation}
for all $m<n_\d$. Since $\a_n \le \a_0$ for all $n$, it follows from (\ref{10.6.2}) that $n_\d<\infty$. The proof is therefore complete.
}

\begin{remark}
{\rm By taking $m =n_\d-1$ in (\ref{10.6.2}), the integer $n_\d$ defined by Rule \ref{Rule1} can be estimated by
\begin{equation}\label{10.6.1}
c_3 \d^2 \sum_{n=0}^{n_\d-1} \frac{1}{\a_n} \le D_{\xi_0} \Theta(\hat x, x_0),
\end{equation}
where $c_3 := c_2 \min\{\mu_0, \mu_1\} \tau^2$. In case $\{\a_n\}$ is chosen such that $\a_{n+1}/\a_n \le q$ for all $n$ for some constant $0<q<1$, then
\begin{equation*}
\sum_{n=0}^{n_\d-1} \frac{1}{\a_n}\ge \frac{1}{\a_0} \sum_{n=0}^{n_\d-1} q^{-n} = \frac{1-q^{n_\d}}{\a_0(1-q)q^{n_\d-1}} \ge \frac{1}{\a_0 q^{n_\d-1}}.
\end{equation*}
It then follows from (\ref{10.6.1}) that $c_3 \a_0^{-1} \d^2 q^{-n_\d+1} \le D_{\xi_0}\Theta(\hat x, x_0)$ which implies that
$n_\d = O(1+|\log \d|)$. Therefore, with such a chosen $\{\a_n\}$, Algorithm \ref{alg1} exhibits the fast convergence property.
}
\end{remark}

In order to use the results given in Lemma \ref{lem1} and Theorem \ref{thm1} to prove the convergence of the method (\ref{method_noise}),
we need the following stability result.

\begin{lemma}\label{lem10.9.1}
Let $\{\xi_n\}$ and $\{x_n\}$ be defined by (\ref{method}) with $\{t_n\}$ chosen by (\ref{10.13.7}), and let $\{\xi_n^\d\}$ and
$\{x_n^\d\}$ be defined by (\ref{method_noise}) with $\{t_n^\d\}$ chosen by (\ref{10.3.6}). Then for each fixed integer $n$ there hold
\begin{equation*}
\lim_{\d \rightarrow 0} \|x_n^\d-x_n\| =0 \quad \mbox{ and } \quad \lim_{\d\rightarrow 0} \|\xi_n^\d -\xi_n\| =0.
\end{equation*}
\end{lemma}

\proof{
We prove the result by induction on $n$. It is trivial when $n=0$ because $\xi_0^\d =\xi_0$ and $x_0^\d =x_0$. Assume next that
the result is true for some $n\ge 0$. We will show that $\xi_{n+1}^\d \rightarrow \xi_{n+1}$ and $x_{n+1}^\d \rightarrow x_{n+1}$
as $\d \rightarrow 0$. We consider two cases:

{\it Case 1: $A x_n \ne y$.}  In this case we must have $A^*(\a_n I +AA^*)^{-1} (A x_n-y)\ne 0$ since otherwise
\begin{equation*}
\fl 0=\l A^*(\a_n I +AA^*)^{-1} (A x_n-y), x_n-x^\dag\r =\|(\a_n I +AA^*)^{-1/2} (A x_n-y)\|^2>0.
\end{equation*}
Therefore, by the induction hypothesis it is straightforward to see that $t_n^\d\rightarrow t_n$ as $\d\rightarrow 0$.
According to the definition of $\xi_{n+1}^\d$ and the induction hypothesis, we then obtain $\lim_{\d\rightarrow 0} \|\xi_{n+1}^\d-\xi_{n+1}\|=0$.
Recall that
\begin{equation*}
x_{n+1}=\nabla \Theta^*(\xi_{n+1}) \quad \mbox{ and } \quad x_{n+1}^\d =\nabla \Theta^*(\xi_{n+1}^\d).
\end{equation*}
We then obtain $\lim_{\d\rightarrow 0} \|x_{n+1}^\d-x_{n+1}\|=0$ by the continuity of $\nabla \Theta^*$.

{\it Case 2: $A x_n=y$.} In this case we have $\xi_{n+1}=\xi_n$. Therefore
\begin{equation*}
\xi_{n+1}^\d-\xi_{n+1} =\xi_n^\d -\xi_n - t_n^\d A^* (\a_n I +AA^*)^{-1} (A x_n^\d-y^\d).
\end{equation*}
Consequently, by the induction hypothesis, we have
\begin{eqnarray*}
\limsup_{\d\rightarrow 0} \|\xi_{n+1}^\d-\xi_{n+1}\|
&\le \limsup_{\d\rightarrow 0} \left(\|\xi_n^\d -\xi_n\|+ \frac{\mu_1}{\sqrt{\a_n}} \|A x_n^\d-y^\d\|\right)\\
& = \frac{\mu_1}{\sqrt{\a_n}} \|A x_n-y\| =0.
\end{eqnarray*}
By using again the continuity of $\nabla \Theta^*$, we obtain $\lim_{\d\rightarrow 0} \|x_{n+1}^\d-x_{n+1}\|=0$.
}

We are now in a position to give the main result concerning the regularization property of the method (\ref{method_noise}) with noisy data
when it is terminated by Rule \ref{Rule1}.

\begin{theorem}\label{thm2}
Let $\Theta: \X \to (-\infty, \infty]$ be proper, lower semi-continuous and strong convex in the sense of (\ref{sc}).
Let $\{\a_n\}$ be a decreasing sequence of positive numbers and let $\{t_n^\d\}$ be chosen by (\ref{10.3.6}) with $0<\mu_0<4 c_0 (1-1/\tau)$
and $\mu_1>0$. Let $n_\d$ be the finite integer defined by Rule \ref{Rule1}. Then for the method (\ref{method_noise})
there hold
\begin{equation*}
\lim_{\d\rightarrow 0} \|x_{n_\d}^\d-x^\dag\| =0 \qquad \mbox{and}
\qquad \lim_{\d\rightarrow 0} D_{\xi_{n_\d}^\d}\Theta (x^\dag, x_{n_\d}^\d) =0.
\end{equation*}
\end{theorem}

\proof{
Due to the strong convexity of $\Theta$, it suffices to show that $\lim_{\d\rightarrow 0} D_{\xi_{n_\d}^\d} \Theta(x^\dag, x_{n_\d}^\d)=0$.
By the subsequence-subsequence argument, we may complete the proof by considering two cases.

Assume first that $\{y^{\d_k}\}$ is a sequence satisfying $\|y^{\d_k}-y\|\le \d_k$ with $\d_k\rightarrow 0$ such that
$n_k:=n_{\d_k}\rightarrow \hat n$ as $k\rightarrow \infty$ for some finite integer $\hat n$. We may assume $n_k=\hat n$ for all $k$.
From the definition of $\hat n:=n_k$ we have
\begin{equation}\label{eq101}
\sqrt{\a_{\hat n}} \|(\a_{\hat n} I +AA^*)^{-1/2} (A x_{\hat n}^{\d_k})-y^{\d_k})\|\le \tau \d_k.
\end{equation}
By taking $k\rightarrow \infty$ and using Lemma \ref{lem10.9.1}, we can obtain $A x_{\hat n}=y$. In view of the definition of $\{\xi_n\}$
and $\{x_n\}$, this implies that $\xi_n=\xi_{\hat n}$ and $x_n=x_{\hat n}$ for all $n\ge \hat n$. Since Theorem \ref{thm1} implies that
$x_n\rightarrow x^\dag$ as $n\rightarrow \infty$, we must have $x_{\hat n}=x^\dag$. Moreover, by Lemma \ref{lem10.9.1},
$\xi_{n_k}^{\d_k} \rightarrow \xi_{\hat n}$ as $k\rightarrow \infty$. Therefore, by the continuity of $\Theta^*$ we can obtain
\begin{eqnarray*}
\lim_{k\rightarrow \infty} D_{\xi_{n_k}^{\d_k}}\Theta(x^\dag, x_{n_k}^{\d_k})
& = \lim_{k\rightarrow \infty} \left(\Theta(x^\dag) + \Theta^*(\xi_{n_k}^{\d_k})- \l \xi_{n_k}^{\d_k}, x^\dag\r\right)\\
& = \Theta(x_{\hat n})+\Theta^*(\xi_{\hat n})-\l \xi_{\hat n}, x_{\hat n}\r=0.
\end{eqnarray*}

Assume next that $\{y^{\d_k}\}$ is a sequence satisfying $\|y^{\d_k}-y\|\le \d_k$ with $\d_k\rightarrow 0$ such that
$n_k:=n_{\d_k}\rightarrow \infty$ as $k\rightarrow \infty$. Let $n$ be any fixed integer. Then $n_k>n$ for large $k$.
It then follows from (\ref{monotone}) in Lemma \ref{lem1} that
\begin{eqnarray*}
D_{\xi_{n_k}^{\d_k}} \Theta(x^\dag, x_{n_k}^{\d_k}) \le D_{\xi_n^{\d_k}} \Theta (x^\dag, x_n^{\d_k})
=\Theta(x^\dag) + \Theta^*(\xi_n^{\d_k}) -\l \xi_n^{\d_k}, x^\dag \r.
\end{eqnarray*}
By using Lemma \ref{lem10.9.1} and the continuity of $\Theta^*$ we obtain
\begin{eqnarray*}
\limsup_{k\rightarrow \infty} D_{\xi_{n_k}^{\d_k}} \Theta (x^\dag, x_{n_k}^{\d_k})
& \le \Theta(x^\dag)+ \Theta^*(\xi_n) - \l \xi_n, x^\dag\r = D_{\xi_n}\Theta(x^\dag, x_n).
\end{eqnarray*}
Since $n$ can be arbitrary and since Theorem \ref{thm1} implies that $D_{\xi_n}\Theta(x^\dag, x_n)\rightarrow 0$
as $n\rightarrow \infty$, we therefore have $\lim_{k\rightarrow \infty} D_{\xi_{n_k}^{\d_k}} \Theta (x^\dag, x_{n_k}^{\d_k})=0$.
}

\begin{remark}
{\rm
In certain applications, the solution of (\ref{1.1}) may have some physical constraints. Thus, instead of (\ref{1.1}), we need to
consider the constrained problem
\begin{equation*}
A x = y \qquad \mbox{subject to } x\in {\mathcal C},
\end{equation*}
where ${\mathcal C}$ is a closed convex subset in $\X$. Correspondingly, (\ref{method_noise}) can be modified into the form
\begin{equation}\label{10.16.1}
\left\{\begin{array}{lll}
\xi_{n+1}^\d  =\xi_n^\d -t_n^\d A^* (\a_n I +AA^*)^{-1} (A x_n^\d -y^\d),\\
x_{n+1}^\d  = \arg \min_{x\in {\mathcal C}} \left\{\Theta(x) -\l \xi_{n+1}^\d, x\r\right\}
\end{array}\right.
\end{equation}
which can be analyzed by the above framework by introducing $\Theta_{\mathcal C}: = \Theta + \imath_{\mathcal C}$, where
$\imath_{\mathcal C}$ denotes the indicator function of ${\mathcal C}$, i.e.
\begin{equation*}
\imath_{\mathcal C}(x) = \left\{\begin{array}{lll}
0,    & x\in {\mathcal C},\\
+\infty, \quad & x \not \in {\mathcal C}
\end{array}\right.
\end{equation*}
When $t_n^\d$ is chosen by (\ref{10.3.6}) and (\ref{10.16.1}) is terminated by Rule \ref{Rule1}, we still have $\|x_{n_\d}^\d-x^\dag\|\rightarrow 0$
and $D_{\xi_{n_\d}^\d} \Theta_{\mathcal C}(x^\dag, x_{n_\d}^\d) \rightarrow 0$ as $\d\rightarrow 0$. However, $D_{\xi_{n_\d}^\d} \Theta(x^\dag, x_{n_\d}^\d)$
may not converge to $0$ because $\xi_{n_\d}^\d$ is not necessarily in $\p \Theta(x_{n_\d}^\d)$.
}
\end{remark}

\begin{remark}
{\rm
In order to implement Algorithm \ref{alg1}, a key ingredient is to solve the minimization problem
\begin{equation}\label{eq4.16.1}
x=\arg \min_{z\in \X} \left\{\Theta(z) -\l \xi, z\r\right\}
\end{equation}
for any given $\xi\in \X$. For some choices of $\Theta$, this minimization problem can be efficiently solved numerically.
When $\X=L^2(\Omega)$, where $\Omega$ is a bounded Lipschitz domain in Euclidean space, there are at least two important choices of $\Theta$
that are crucial for sparsity recovery and discontinuity detection. The first one is
\begin{equation}\label{eq:L1}
\Theta(x):= \frac{1}{2\beta} \int_\Omega |x(\omega)|^2 d\omega +\int_\Omega |x(\omega| d\omega
\end{equation}
with $\beta>0$, the minimizer of (\ref{eq4.16.1}) can be given explicitly by
\begin{equation*}
x(\omega)=\beta \mbox{sign}(\xi(\omega)) \max\{ |\xi(\omega)|-1, 0\}, \qquad \omega \in \Omega.
\end{equation*}
The second one is
\begin{equation}\label{eq:TV}
\Theta(x):= \frac{1}{2\beta} \int_\Omega |x(\omega)|^2 d\omega +\mbox{TV}(x)
\end{equation}
with $\beta>0$, where $\mbox{TV}(x)$ denotes the total variation of $x$, i.e.
\begin{equation*}
\mbox{TV}(x) :=\sup\left\{ \int_\Omega x \, \mbox{div} f d\omega:
f\in C_0^1(\Omega; {\mathbb R}^N) \mbox{ and }
\|f\|_{L^\infty(\Omega)}\le 1\right\}.
\end{equation*}
Then the minimization problem (\ref{eq4.16.1}) can be equivalently formulated as
\begin{equation*}
x =\arg \min_{z\in L^2(\Omega)}  \left\{ \frac{1}{2\beta} \|z-\beta \xi\|^2_{L^2(\Omega)} + \mbox{TV}(z) \right\}
\end{equation*}
which is the total variation denoising problem (\cite{ROF92}). Although there is no explicit formula for the minimizer of (\ref{eq4.16.1}),
there are many efficient numerical solvers developed in recent years, see \cite{BT2009a,BT2009b,CP2011,MSX2011}. For the numerical simulations
involving total variation presented in Section 3, we always use the denoising algorithm \texttt{FISTA} from \cite{BT2009a,BT2009b}.
Indeed, when solving (\ref{eq4.16.1}) with $\Theta$ given by (\ref{eq:TV}), \texttt{FISTA} is used to solve its dual problem whose solution
determines the solution of the primal problem (\ref{eq4.16.1}) directly; one may refer to the algorithm on page 2462 in \cite{BT2009b} and its
monotone version\begin{footnote}
{After acceptance of this paper, we found that our method can be significantly accelerated if we use PDHG (an application of Uzawa algorithm)
to solve the TV denoising problem. The Matlab code of PDHG can be found at http://pages.cs.wisc.edu/$\sim$swright/GPUreconstruction/}
\end{footnote}.
}
\end{remark}

\begin{remark}
{\rm
Another key ingredient in implementing Algorithm \ref{alg1} is to determine $v := (\a I +AA^*)^{-1} r$
for $\a>0$, where $r:=A x_n^\d -y^\d$. This amounts to solving the linear equation
\begin{equation*}
(\a I +AA^*) v =r
\end{equation*}
for which many efficient solvers from numerical linear algebra can be applied. When $A$ has special structure, this equation
can even be solved very fast. For instance, if $A$ is a convolution operator in ${\mathbb R}^d$, say
\begin{equation*}
A x(\sigma) =\int_{{\mathbb R}^d} k(\sigma-\eta) x(\eta) d\eta
\end{equation*}
with the kernel $k$ decaying sufficiently fast at infinity, then $v$ can be determined as
\begin{equation*}
v = {\mathcal F}^{-1} \left(\frac{{\mathcal F}(r)}{\a +|{\mathcal F}(k)|^2}\right),
\end{equation*}
where ${\mathcal F}$ and ${\mathcal F}^{-1}$ denote the Fourier transform and the inverse Fourier transform respectively.
Therefore $v$ can be calculated efficiently by the fast Fourier transform.
}
\end{remark}

\subsection{\bf Possible extension for nonlinear inverse problems}

Our method can be extended for solving nonlinear inverse problems in Hilbert spaces that can be formulated as the equation
\begin{equation}\label{n1.1}
F(x) =y,
\end{equation}
where $F: \mathscr{D}(F) \subset \X \to \Y$ is a nonlinear continuous operator between two Hilbert spaces $\X$ and $\Y$ with closed convex
domain $\mathscr{D}(F)$. We assume that for each $x \in \mathscr{D}(F)$ there is a bounded linear operator $L(x): \X \to \Y$ such that
\begin{equation*}
\lim_{h\searrow 0} \frac{F(x+h(z-x))-F(x)}{h} = L(x) (z-x), \quad \forall z\in \mathscr{D}(F).
\end{equation*}
In case $F$ is Fr\'{e}chet differentiable at $x$, $L(x)$ is exactly the Fr\'{e}chet derivative of $F$ at that point.

In order to find the solution of (\ref{n1.1}) with special feature, as before we introduce a penalty function $\Theta: \X \to (-\infty, \infty]$
which is proper, convex and lower semi-continuous.  Let $y^\d$ be the only available noisy data satisfying
\begin{equation*}
\|y^\d-y\|\le \d
\end{equation*}
with a small known noise level $\d>0$. Then an obvious extension of Algorithm \ref{alg1} for solving (\ref{n1.1}) takes the
following form.

\begin{algorithm}[Nonstationary iterative method for nonlinear problem] \label{alg2} \quad
\begin{enumerate}
\item[(i)] Take $\tau>1$, $\mu_0>0$, $\mu_1>0$ and a decreasing sequence $\{\a_n\}$ of positive numbers;

\item[(ii)] Take $\xi_0\in \X$ and define $x_0:=\arg\min_{x\in D(F)} \{\Theta(x) -\l \xi_0, x\r\}$ as an initial guess;

\item[(iii)]  For each $n=0, 1, \cdots$ define
\begin{eqnarray*}
\xi_{n+1} &= \xi_n - t_n L(x_n)^* (\a_n I +L(x_n) L(x_n)^*)^{-1} (F(x_n) -y^\d),\\
x_{n+1} & = \arg \min_{x\in D(F)} \{ \Theta(x) -\l \xi_{n+1}, x\r\},
\end{eqnarray*}
where
\begin{equation*}
\fl \qquad t_n = \min \left\{ \frac{\mu_0 \l (\a_n I + L(x_n)L(x_n)^*)^{-1} (F(x_n)-y^\d), F(x_n)-y^\d\r}
{\|L(x_n)^*(\a_n I +L(x_n) L(x_n)^*)^{-1} (F(x_n)-y^\d)\|^2}, \mu_1\right\};
\end{equation*}

\item[(iv)] Let $n_\d$ be the first integer such that
\begin{equation*}
\fl \qquad \a_{n_\d} \l (\a_{n_\d} I + L(x_{n_\d}) L(x_{n_\d})^*)^{-1} (F(x_{n_\d}) -y^\d), F(x_{n_\d})-y^\d\r \le \tau^2 \d^2
\end{equation*}
and use $x_{n_\d}$ to approximate the solution of (\ref{n1.1}).
\end{enumerate}
\end{algorithm}

We remark that when $\Theta(x) =\frac{1}{2} \|x\|^2$, Algorithm \ref{alg2} reduces to a method which is similar to the regularized
Levenberg-Marquardt method in \cite{Jin2009} for which convergence is proved under certain conditions on $F$. For general convex
penalty function $\Theta$, however, we do not have convergence theory on Algorithm \ref{alg2} yet. Nevertheless, we will give numerical
simulations to indicate that it indeed performs very well.

\section{Numerical simulations}

In this section we will provide various numerical simulations on our method. The choice of the sequence $\{\a_n\}$ plays a crucial role
for the performance: if $\{\a_n\}$ decays faster, only fewer iterations are required but the reconstruction result is less accurate; on the other hand,
if $\{\a_n\}$ decays slower, more iterations are required but the reconstruction result is more accurate. In order to solve this dilemma,
we choose fast decaying $\{\a_n\}$ at the beginning, and then choose slow decaying $\{\a_n\}$ when the method tends to stop. More precisely, we
choose $\{\a_n\}$ according to the following rule.

\begin{Rule}\label{alpha}
Let $0<\gamma_0\le \gamma_1\le 1$ and $\hat \rho>1$ be preassigned numbers. We take some number $\a_0>0$ and for $n\ge 0$ define
\begin{equation*}
\rho_n := \frac{\sqrt{\a_n} \| (\a_n I +AA^*)^{-1/2}(A x_n^\d-y^\d)\|}{\tau\d}.
\end{equation*}
If $\rho_n > \hat \rho$ we set $\a_{n+1} =\gamma_0\a_n$; otherwise we set $\a_{n+1}=\gamma_1\a_n$.
\end{Rule}

All the computation results in this section are based on $\{\a_n\}$ chosen by this rule with $\gamma_0 \approx 0.5$, $\gamma_1 \approx 1$
 and $\hat \rho \approx 2.5$. Our tests were done by using MATLAB R2012a on an Lenovo laptop with Intel(R)
Core(TM) i5 CPU 2.30 GHz and 6 GB memory.

\subsection{\bf Integral equation of first kind in dimension one}

We first consider the integral equation of the form
\begin{equation}\label{Int}
Ax(s):=\int_0^1 k(s,t) x(t) dt =y(s) \quad \mbox{ on } [0,1],
\end{equation}
where
\begin{equation*}
k(s,t)=\left\{\begin{array}{lll}
40 s(1-t), & s\le t\\
40 t(1-s), & s\ge t.
\end{array}\right.
\end{equation*}
It is easy to see, that $A$ is a compact linear operator from $L^2[0,1]$ to $L^2[0,1]$.
Our goal is to find the solution of (\ref{Int}) using noisy data $y^\d$ satisfying $\|y-y^\d\|_{L^2[0,1]}=\d$ for some
specified noise level $\d$. In our numerical simulations, we divide $[0,1]$ into $N=400$ subintervals of equal length
and approximate any integrals by the trapezoidal rule.

\begin{figure}[ht]
\centering
  \includegraphics[width = 1\textwidth, height=3in]{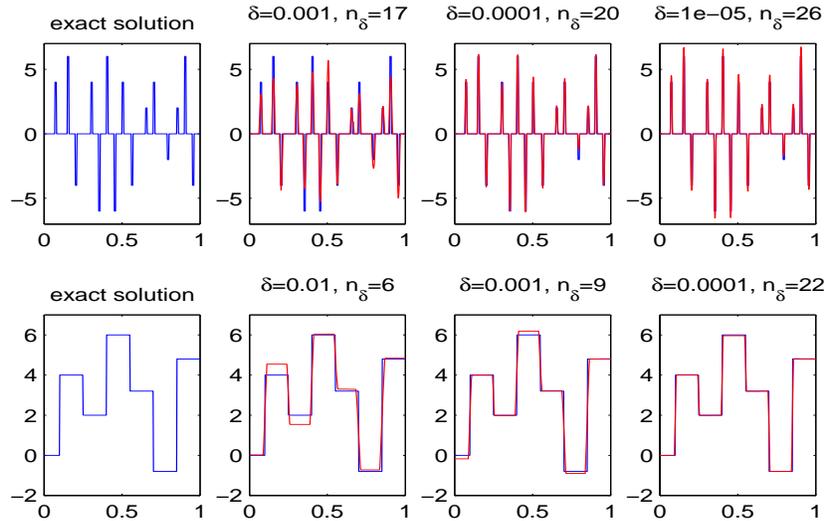}
  \caption{Reconstruction results for (\ref{Int}) by our method using noisy data with various noise levels}\label{fig1}
\end{figure}

In Figure \ref{fig1} we report the numerical performance of Algorithm \ref{alg1}. The sequence $\{\a_n\}$ is selected by Rule \ref{alpha} with
$\a_0=0.01$, $\gamma_0=0.6$, $\gamma_1 =0.99$ and $\hat \rho =2.5$. The first row gives the reconstruction results using noisy data
with various noise levels when the sought solution is sparse; we use the penalty function $\Theta$ given in (\ref{eq:L1}) with $\beta =10$.
The second row reports the reconstruction results for various noise levels when the sought solution is piecewise constant; we use the
penalty function $\Theta$ given in (\ref{eq:TV}) with $\beta=100$. When the 1d TV-denoising algorithm \texttt{FISTA}
in \cite{BT2009a,BT2009b} is used to solve the minimization problems associated with this $\Theta$, it is terminated as long as the number of iterations
exceeds $2500$ or the error between two successive iterates is smaller than $10^{-6}$.
During these computations, we use $\xi_0(t) \equiv 0$ and the parameters $\tau =1.01$, $\mu_0 = 1/\beta$ and $\mu_1 =1$ in Algorithm \ref{alg1}.
The computational times for the first row are $0.0677$, $0.0826$ and $0.1017$ seconds respectively, and the computation times for the second row
are $0.2963$, $0.4458$ and $1.0672$ seconds respectively.
This shows that Algorithm \ref{alg1} indeed is a fast method with the capability of capturing special features of solutions.

\subsection{\bf Determine source term in Poisson equation}

Let $\Omega = [0,1]\times [0,1]$. We consider the problem of determining the source term $f\in L^2(\Omega)$ in the Poisson equation
\begin{equation*}
-\triangle u = f \quad \mbox{in } \Omega, \qquad   u =0 \quad \mbox{on } \p \Omega
\end{equation*}
from an $L^2(\Omega)$ measurement $u^\d$ of $u$ with $\|u^\d-u\|_{L^2(\Omega)} \le \d$. This problem takes the form (\ref{1.1}) if
we define $A:=(-\triangle)^{-1}$, where $-\triangle: H^2\cap H_0^1(\Omega) \to L^2(\Omega)$ is an isomorphism. The information on
$A$ can be obtained by solving the equation.

In order to solve the Poisson equation numerically, we take $(N+1)\times (N+1)$ grid points
\begin{equation*}
(x_i, y_j) := (i/N, j/N), \qquad i, j=0, 1, \cdots, N
\end{equation*}
on $\Omega$, and write ${\bf u}_{i,j}$ for $u(x_i, y_j)$ and ${\bf f}_{i,j}$ for $f(x_i, y_j)$.
By the finite difference representation of $-\triangle u$, the Poisson equation has the discrete form
\begin{equation}\label{10.29.1}
\fl \qquad 4 {\bf u}_{i,j}-{\bf u}_{i+1, j} -{\bf u}_{i-1,j} -{\bf u}_{i,j+1}-{\bf u}_{i,j-1} = h^2 {\bf f}_{i,j}, \quad i, j=1, \cdots, N-1,
\end{equation}
where $h=1/N$. Since $u=0$ on $\p \Omega$, the discrete sine transform can be used to solve (\ref{10.29.1}). Consequently ${\bf u}_{i,j}$
can be determined by the inverse discrete sine transform (\cite{PTVF2007})
\begin{equation*}
{\bf u}_{i,j} = ({\bf S}^{-1} \hat{\bf u})_{i,j} :=4 h^2 \sum_{p=1}^{N-1} \sum_{q=1}^{N-1} \hat {\bf u}_{p,q} \sin(iph\pi) \sin(jqh\pi)
\end{equation*}
for $i, j=1, \cdots, N-1$, where
\begin{equation*}
\hat {\bf u}_{p,q} = ({\bf \Lambda} \hat {\bf f})_{p,q}:= \frac{h^2 \hat {\bf f}_{p,q}}{4-2\cos(ph\pi)-2\cos(qh\pi)}
\end{equation*}
and $\hat {\bf f}_{p,q}$ is determined by the discrete sine transform
\begin{equation*}
\hat {\bf f}_{p,q} = ({\bf S} {\bf f})_{p,q}:= \sum_{i=1}^{N-1} \sum_{j=1}^{N-1} {\bf f}_{i,j} \sin(iph\pi) \sin(jqh\pi).
\end{equation*}
Let ${\bf A} = {\bf S}^{-1} {\bf \Lambda} {\bf S}$. Then ${\bf f}$ can be determined by solving the equation ${\bf A} {\bf f} = {\bf u}$.
When applying Algorithm \ref{alg1}, we need to determine ${\bf v}= (\a I +{\bf A}{\bf A}^*)^{-1} {\bf r}$ for various $\a>0$ and vectors ${\bf r}$.
This can be computed as
\begin{equation*}
{\bf v} = {\bf S} (\a I +{\bf \Lambda}^2)^{-1} {\bf S}^{-1} {\bf r},
\end{equation*}
where, for any vector ${\bf w}$,  ${\bf S} {\bf w}$ and ${\bf S}^{-1} {\bf w}$ can be implemented by the fast sine and inverse sine transforms respectively, while
\begin{equation*}
[(\a I +{\bf \Lambda}^2)^{-1} {\bf w}]_{i,j}= {\bf w}_{i,j}/(\a +h^4/(4-2\cos(ih\pi)-2\cos(jh\pi))^2)
\end{equation*}
Therefore ${\bf v}$ can be computed efficiently.

\begin{figure}[ht]
\centering
  \includegraphics[width = 1\textwidth, height=2in]{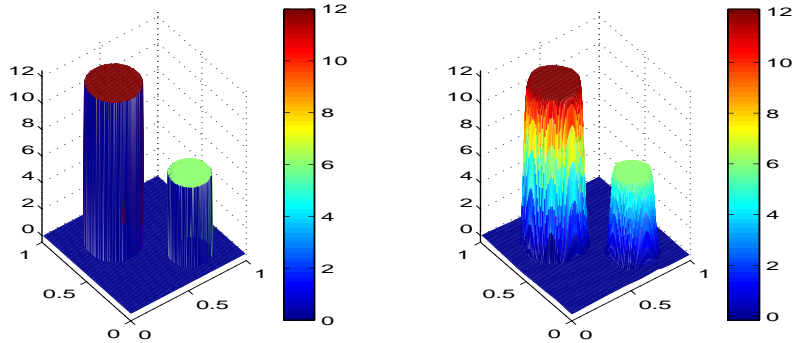}
  \caption{Reconstruction of the source term in Poisson equation using noisy data with $\delta =10^{-3}$}\label{fig4}
\end{figure}

We apply Algorithm \ref{alg1} to reconstruct the source term which is assumed to be piecewise constant. In our computation we use a noisy
data with noise level $\d=10^{-3}$. The left plot in Figure \ref{fig4} is the exact solution. The right plot in Figure \ref{fig4} is the
reconstruction result by Algorithm \ref{alg1} using initial guess $\xi_0 \equiv 0$ and the penalty function
\begin{equation*}
\Theta({\bf f}) = \frac{1}{2\beta} \|{\bf f}\|_F + \mbox{TV}_I({\bf f}),
\end{equation*}
where $\|{\bf f}\|_F$ is the Frobenius norm of ${\bf f}$ and $\mbox{TV}_I({\bf f})$ denotes the discrete isotropic TV defined by  (\cite{BT2009b})
\begin{eqnarray*}
\mbox{TV}_I({\bf f}) & := \sum_{i=1}^{N-1} \sum_{j=1}^{N-1} \sqrt{({\bf f}_{i,,j}-{\bf f}_{i+1,j})^2+({\bf f}_{i,j}-{\bf f}_{i,j+1})^2} \\
& \quad + \sum_{i=1}^{N-1} |{\bf f}_{i, N} -{\bf f}_{i+1,N}| +\sum_{j=1}^{N-1} |{\bf f}_{N,j}-{\bf f}_{N, j+1}|.
\end{eqnarray*}
In each step of Algorithm \ref{alg1}, the minimization problem associate with $\Theta$ is solved by performing 400 iterations of the 2d TV-denoising
algorithm \texttt{FISTA} in \cite{BT2009b}. In our computation, we use $N=120$, and for those parameters in Algorithm \ref{alg1}, we take $\tau =1.01$, $\beta = 20$,
$\mu_0 = 0.4/\beta$ and $\mu_1 =2$. When using Rule \ref{alpha} to choose $\{\a_n\}$ we take $\a_0 = 0.001$, $\gamma_0=0.5$,
$\gamma_1 = 0.95$ and $\hat \rho = 2$. The reconstruction result indicates that our method succeeds in capturing the feature of the solution. Moreover, the computation terminates
after $n_\d=17$ iterations and takes 11.7092 seconds.

\subsection{\bf Image deblurring}

Blurring in images can arise from many sources, such as limitations of the optical system, camera and object motion, astigmatism,
and environmental effects (\cite{HNO2006}). Image deblurring is the process of making a blurry image clearer to better represent the
true scene.

We consider grayscale digital images which can be represented as rectangular matrices of size $m \times n$.
Let ${\bf X}$ and ${\bf B}$ denote the true image and the blurred image respectively. The blurring process can be
described by an operator ${\bf L}: {\mathbb R}^{m\times n}\to {\mathbb R}^{m\times n}$ such that ${\bf B} ={\bf L} ({\bf X})$.
We consider the case that the model is shift-invariant and ${\bf L}$ is linear. By stacking the columns of ${\bf X}$ and ${\bf B}$
we can get two long column vectors ${\bf x}$ and ${\bf b}$ of length $N:=mn$. Then there is a large matrix
${\bf A} \in {\mathbb R}^{N\times N}$ such that ${\bf b} = {\bf A} {\bf x}$. Considering the appearance of unavoidable random noise,
one in fact has
\begin{equation*}
{\bf b}^\d = {\bf A} {\bf x} + {\bf e},
\end{equation*}
where ${\bf e}$ denotes the noise. The blurring matrix ${\bf A}$ is determined by the point spread function (PSF) ${\bf P}$---the function that
describes the blurring and the resulting image of the single bright pixel (i.e. point source).

Throughout this subsection, periodic boundary conditions are assumed on all images. Then ${\bf A}$ is a matrix which is block circulant with circulant blocks;
each block is built from ${\bf P}$. It turns out that ${\bf A}$ has the spectral decomposition
\begin{equation*}
{\bf A} ={\bf F}^* {\bf \Lambda} {\bf F}
\end{equation*}
where {\bf F} is the two-dimensional unitary discrete Fourier transform matrix and ${\bf \Lambda}$ is the diagonal matrix
whose diagonal entries are eigenvalues of ${\bf A}$. The diagonal matrix ${\bf \Lambda}$ is easily determined by the smaller matrix ${\bf P}$,
and the action of ${\bf F}$ and ${\bf F}^*$ can be realized by \texttt{fft} and \texttt{ifft}. Therefore, for any ${\bf v} \in {\mathbb R}^N$
and $\a>0$, $(\a {\bf I} + {\bf A}{\bf A}^*)^{-1} {\bf v}$ is easily computable by the fast Fourier transform.

\begin{figure}[ht]
\centering
  \includegraphics[width = 1\textwidth, height=3.6in]{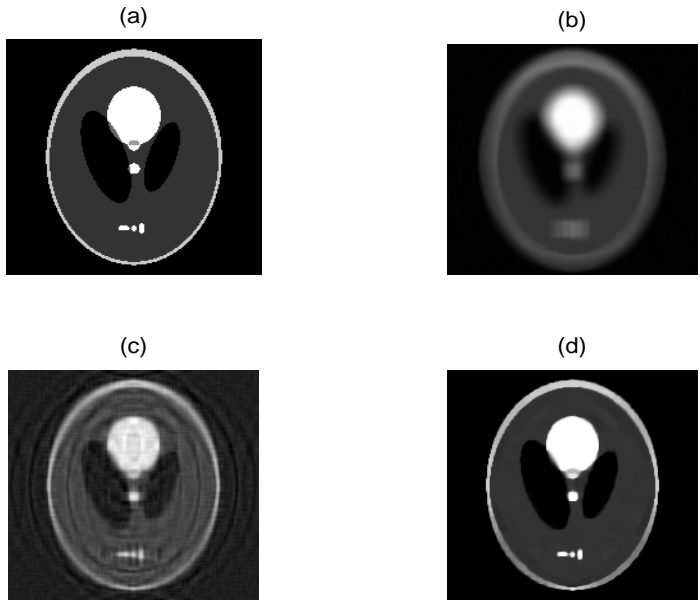}
  \caption{Reconstruction of the Shepp-Logan phantom of size $200\times 200$ blurred by Gaussian PSF: (a) original image;
  (b) blurred image corrupted by Gaussian noise with $\delta_{rel} =1.25\%$; (c) reconstruction result by Algorithm \ref{alg1}
  with $\Theta=\Theta_0$; (d) reconstruction result by Algorithm \ref{alg1} with $\Theta =\Theta_1$. }\label{fig5}
\end{figure}

In the following we perform some numerical experiments by applying Algorithm \ref{alg1} to deblur various corrupted images. In our simulations
the exact data ${\bf b}$ are contaminated by random noise vectors ${\bf e}$ whose entries are normally distributed with zero mean.
We use
$$
\d_{rel} := \frac{\|{\bf e}\|_2}{\|{\bf b}\|_2}
$$
to denote the relative noise level. When applying Algorithm \ref{alg1}, we use $\xi_0=0$ and the following two convex functions
$$
\Theta_0({\bf X}) =\frac{1}{2} \|{\bf X}\|_F^2 \qquad \mbox{and} \qquad \Theta_1({\bf X}) =\frac{1}{2} \|{\bf X}\|_F^2 + \mbox{TV}_I({\bf X}),
$$
For those parameters in the algorithm we take $\tau =1.001$, $\mu_0 =0.4$ and $\mu_1 =2$. In each step of the algorithm, the minimization problem
associate with $\Theta$ is solved by performing 200 iterations of the algorithm \texttt{FISTA}
in \cite{BT2009b}. When using Rule \ref{alpha} to choose $\{\a_n\}$ we take $\a_0 = 1$, $\gamma_0=0.5$,
$\gamma_1 = 0.99$ and $\hat \rho = 2.5$. In order to compare the quality of the restoration $\widetilde{\bf X}$, we evaluate the
peak signal-to-noise ratio (PSNR) value defined by
$$
PSNR = 20 \log_{10} \frac{\sqrt{mn} \max({\bf X})}{\|{\bf X}-\widetilde{\bf X}\|_F},
$$
where $\max({\bf X})$ denotes the maximum possible pixel value of the true image ${\bf X}$.

\begin{figure}[ht]
\centering
  \includegraphics[width = 1\textwidth, height=3.6in]{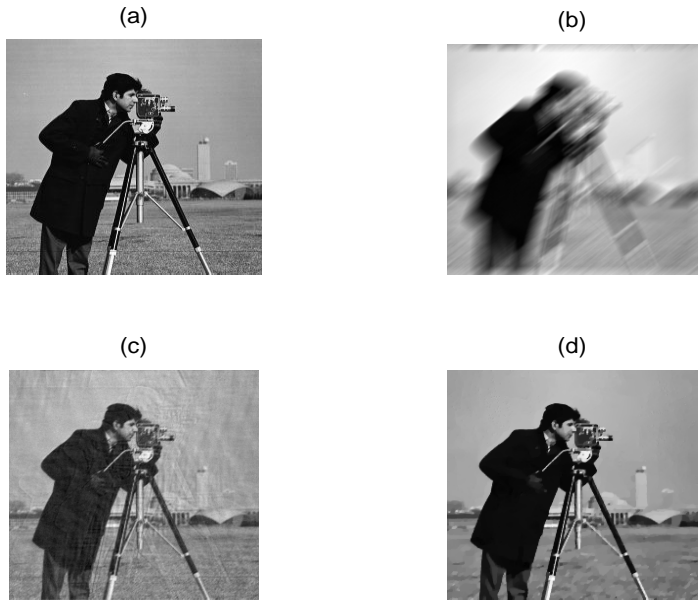}
  \caption{Restoration of the $256\times 256$ Cameraman image blurred by motion: (a) original images; (b) images blurred by motion and noise;
  (c) restoration by Algorithm \ref{alg1} with $\Theta=\Theta_0$;
  (d) restoration by Algorithm \ref{alg1} with $\Theta=\Theta_1$.
  }\label{fig6}
\end{figure}

In Figure \ref{fig5} we plot the restoration results of the Shepp-Logan phantom of size $200\times 200$ which is blurred
by a $15\times 15$ Gaussian PSF with standard derivation $30$ and is contaminated by Gaussian white noise with relative
noise level $\delta_{rel}=1.25\%$. The original and blurred images are plotted in (a) and (b) of Figure \ref{fig5} respectively.
In Figure \ref{fig5} (c) we plot the restoration result by Algorithm \ref{alg1} with $\Theta = \Theta_0$. With such chosen
$\Theta$, the method in Algorithm \ref{alg1} reduces to the classical nonstationary iterated Tikhonov regularization (\ref{NSIT})
which has the tendency to over-smooth solutions. The plot clearly indicates this drawback because of the appearance of the ringing artifacts.
The corresponding PSNR value is $21.3485$. In Figure \ref{fig5} (d) we plot the restoration result by Algorithm \ref{alg1} with $\Theta = \Theta_1$.
Due to the appearance of the total variation term in $\Theta_1$, the artifacts are significantly removed. In fact the corresponding PSNR
value is $24.8653$; the computation terminates after $n_\d=45$ iterations and takes $55.8815$ seconds.

In Figure \ref{fig6} we plot the restoration results of the $256\times 256$ Cameraman image corrupted by a $21\times 25$
linear motion kernel generated by \texttt{fspecial('motion',30,40)} and  a Gaussian white noise with relative
noise level $\delta_{rel}=0.2\%$. The original and blurred images are plotted in (a) and (b) of Figure \ref{fig6} respectively.
In (c) and (d) of Figure \ref{fig6} we plot the restoration results by Algorithm \ref{alg1} with $\Theta = \Theta_0$ and $\Theta =\Theta_1$
respectively. The plot in (c) contains artifacts that degrade the visuality, the plot in (d) removes the artifacts significantly.
In fact the PSNR values corresponding to (c) and (d) are $26.9158$ and $29.8779$ respectively. The computation for (d) terminates
after $n_\d=35$ iterations and takes $76.8809$ seconds.

\subsection{\bf De-autoconvolution}

We finally present some numerical simulations for nonlinear inverse problems by solving the autoconvolution equation
\begin{equation}\label{autoconv}
\int_0^t x(t-s) x(s) ds =y(t)
\end{equation}
defined on the interval $[0, 1]$. The properties of the autoconvolution operator $[F(x)](t): =\int_0^t x(t-s) x(s) ds$ have been discussed in \cite{GH1994}.
In particular, as an operator from $L^2[0,1]$ to $L^2[0,1]$, $F$ is Fr\'{e}chet differentiable; its Fr\'{e}chet derivative and the adjoint
are given respectively by
\begin{eqnarray*}
\,\,\, \,\left[F'(x) v\right](t)  = 2 \int_0^t x(t-s) v(s) ds, \quad v \in L^2[0,1], \\
\left[F'(x)^* w\right](s) = 2 \int_s^1 w(t) x(t-s) dt, \quad w\in L^2[0,1].
\end{eqnarray*}

\begin{figure}[htp]
  \begin{center}
    \includegraphics[width = 1\textwidth, height=3.4in]{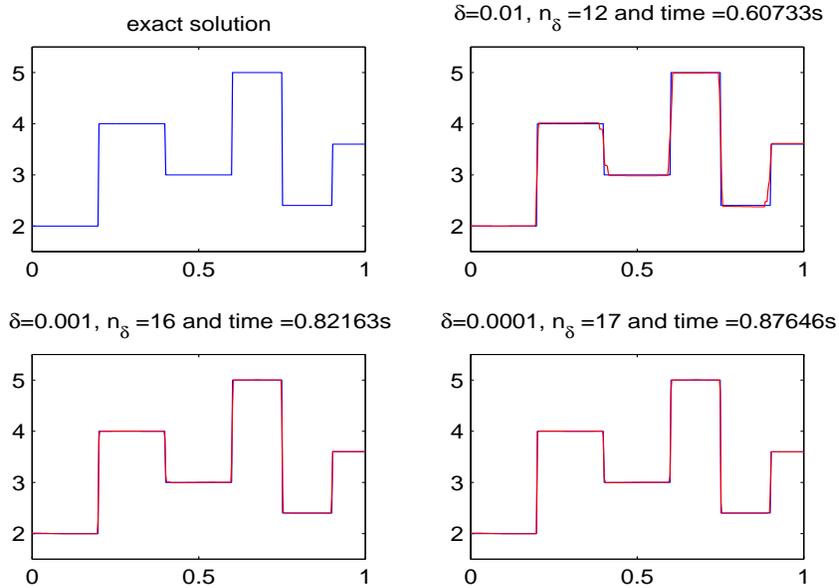}
  \end{center}
  \caption{Reconstruction results for the de-autoconvolution problem by Algorithm \ref{alg2} using noisy data with various noise levels. }\label{fig7}
\end{figure}

We assume that (\ref{autoconv}) has a piecewise constant solution and use a noisy data $y^\d$ satisfying $\|y^\d-y\|_{L^2[0,1]}=\d$ to
reconstruct the solution. In Figure \ref{fig7} we report the reconstruction results by Algorithm \ref{alg2} using $L(x_n)=F'(x_n)$ and
the $\Theta$ given in (\ref{eq:TV}) with $\beta=20$. All integrals involved are approximated by the trapezoidal rule by dividing $[0,1]$
into $N=400$ subintervals of equal length. For those parameters involved in the algorithm, we take $\tau=1.01$, $\mu_0=0.4/\beta$ and $\mu_1=1$.
We also take the constant function $\xi_0(t)\equiv 1/\beta$ as an initial guess. The sequence $\{\a_n\}$ is selected by Rule \ref{alpha} with $A$
replaced by $L(x_n)$ in which $\a_0=1$, $\gamma_0=0.5$, $\gamma_1 =0.99$ and $\hat \rho =3$. When the 1d-denoising algorithm \texttt{FISTA}
in \cite{BT2009a,BT2009b} is used to solve the minimization problems associated with $\Theta$, it is terminated as long as the number of iterations
exceeds $1200$ or the error between two successive iterates is smaller than $10^{-5}$. We indicate in Figure \ref{fig7} the number of iterations
and the computational time for various noise levels $\delta$; the results show that Algorithm \ref{alg2} is indeed a fast method for this problem.

\section{Conclusion}

We proposed a nonstationary iterated method with convex penalty term for solving inverse problems in Hilbert spaces.
The main feature of our method is its splitting character, i.e. each iteration consists of two steps: the first step involves
only the operator from the underlying problem so that the Hilbert space structure can be exploited, while the second step involves
merely the penalty term so that only a relatively simple strong convex optimization problem needs to be solved.
This feature makes the computation much efficient. When the underlying problem is linear, we proved the convergence of our method
in the case of exact data; in case only noisy data are available, we introduced a stopping rule to terminate the iteration and proved
the regularization property of the method. We reported various numerical results which indicate the good performance of our method.
%\begin{figure}[ht]
%\centering
%  \includegraphics[width=5in, height=2.5in]{L1_1d.eps}
%  \caption{Reconstruction results for Example \ref{ex1} with $n_\d$ being the number
%  of iterations. (a) $n_\d=16$; (b) $n_\d=16$; (c) $n_\d=24$.}\label{fig1}
%\end{figure}

\section*{Acknowledgement}

Q. Jin is partially supported by the DECRA grant DE120101707 of Australian Research Council and X. Lu
is partially supported by National Science Foundation of  China (No. 11101316 and No. 91230108).

%%%%%%%%%%%%%%%%%%%%%%%%%%%%%%%%%%%%%%%%%%%%%%%%%%%%%%%%%%%%%%%%%%%%%%%%%%%
%

\end{document}